\documentclass[a4paper]{amsart}

\usepackage{epsfig,amsmath}
\usepackage{xspace}
\usepackage[psamsfonts]{amssymb}
\usepackage[latin1]{inputenc}
\usepackage{float}\restylefloat{figure}
\usepackage{color}

\newtheorem*{lemma*}{\bf Lemma}
\newtheorem*{sublemma*}{\bf Sublemma}
\newtheorem*{claim*}{\bf Claim}
\newtheorem*{complement*}{\bf Complement}

\renewcommand{\epsilon}{\varepsilon}
\renewcommand{\emptyset}{\varnothing}

\newcommand{\cal}[1]{\mathcal{#1}}
\newcommand{\wt}[1]{\widetilde{#1}}
\newcommand{\setof}[2]{\big\{{#1}\,\big|\,{#2}\big\}}

\newcounter{rememberItem}

\def\bEA{\begin{eqnarray*}}
\def\eEA{\end{eqnarray*}}
\def\bEAn{\begin{eqnarray}}
\def\eEAn{\end{eqnarray}}

\def\ds{\displaystyle}
\def\tend{\longrightarrow}

\def\ds{\displaystyle}

\def\on{\operatorname}

\def\cal{\mathcal}

\def\C{{\mathbb C}}
\def\D{{\mathbb D}}
\def\H{{\mathbb H}}
\def\N{{\mathbb N}}

\def\R{{\mathbb R}}

\def\Z{{\mathbb Z}}

\def\Im{{\on{Im}\,}}

\def\remark{\vskip.2cm \noindent{\bf Remark. }}
\def\endremark{\par\medskip}

\newtheorem{theorem}{Theorem}
\newtheorem{lemma}[theorem]{Lemma}

\newtheorem*{theorem*}{Theorem}
\newtheorem*{claim}{Proposal}
\newtheorem*{conjecture}{Conjecture}
\newcommand{\one}{\mathbf{1}}

\newcommand{\Dom}{\mathit{Dom}}
\newcommand{\dom}{\mathit{dom}}
\newcommand{\euh}{E}
\newcommand{\ibar}{{\overline{\imath}}}
\newcommand{\jbar}{{\overline{\jmath}}}

\begin{document}

\title[Siegel disks with non-locally connected boundaries]{Relatively compact Siegel disks with non-locally connected boundaries}
\author{Arnaud Chéritat}

\abstract
We construct holomorphic maps $f$ with a Siegel disk whose boundary is not locally connected (and is an indecomposable continuum), yet compactly contained in the domain of definition of the map. Our examples are injective and defined on a subset of $\C$.
\endabstract

\maketitle

\section{Introduction}

In \cite{handel}, Handel constructed a $C^{\infty}$ area preserving diffeomorphism of the plane that has a minimal set that is a pseudo circle.
In \cite{herman}, Herman adapted the construction to produce a $C^{\infty}$ diffeomorphism of the sphere that is conjugated to a rotation in the two complementary components of an invariant pseudo circle, and holomorphic in one of them. In the same spirit, in \cite{ricardo} Pérez Marco was able, using tube-log Riemann surfaces, to construct examples of injective holomorphic maps defined in a subset $U$ of $\C$ that have a Siegel disk compactly contained in $U$ whose boundary is a $\C^\infty$ Jordan curve, which came as a surprise. Again the method is versatile and Kingshook Biswas used Pérez-Marco's construction to produce a set of interesting examples: \cite{biswas2,biswas1}. Here we add an ingredient to this construction and get:

\begin{theorem} There exists a holomorphic map $f$ defined in a simply connected open subset $U$ of $\C$ containing the origin, fixing $0$ and having at $0$ a Siegel disk $\Delta$ that is compactly contained in $U$ and whose boundary is a pseudo circle.
\end{theorem}

As Herman remarked, the construction is very flexible. The reader will find in the conclusion section more properties that these maps can be given.
There, we also suggest possible other consequences of the method.

\remark We ignore if the maps can be chosen entire (or entire meromorphic).\endremark

\section{The new ingredient}

\subsection{The tool}

Compared to the previous constructions, the main new idea is to make use of Runge's theorem:

\begin{theorem*}[Runge, simply connected case] For all holomorphic map $f$ defined in a simply connected open set $U\subset\C$, there exists a sequence of polynomials that tends to $f$ locally uniformly on $U$.\footnote{a.k.a. ``on every compact subset'' of $U$}
\end{theorem*}

It is easy to see that $f$ can be approximated (in the sense above) by polynomials if and only if it can be approximated by entire maps. In our application, the fact that the approximating sequence consists in polynomials is not so important.
We will only need entire maps.

\subsection{Personal comments}

A while ago, Kingshook Biswas and the author had discussions about the difficulties of the construction of non locally connected relatively compact Siegel disks, and after several failures the latter was intimately convinced of the impossibility to construct such an example. Biswas and Pérez Marco did not look so convinced.

Later, concerning unrelated work, Yohann Genzmer asked the author if he knew how to prove or disprove a result about the radius of convergence of power series with a holomorphic parameter that would have helped him. The author eventually came up with a counterexample, built using iterated exponential maps. Then we asked Julien Duval who explained us how Runge's theorem could more easily yield a counterexample.

More recently, the author decided to read Herman's article on the construction of a smooth diffeomorphism of the sphere, leaving invariant a pseudo circle, conjugated to rotations outside, and holomorphic on one complementary component. The author wondered if Herman's example could be made holomorphic beyond the boundary. While trying again using Pérez Marco's method, he came up using iterated exponential maps. This reminded him of the Genzmer episode. And suddenly it was clear that Runge's theorem was the solution.

It must be acknowledged that Biswas's and Pérez Marco's intuitions were right, after all.

\section{Terminology}

Saying that two holomorphic maps commute can mean different things. 
To be rigorous, we will say that ``$f$ and $g$ commute, domain included'', if the domain of definition of $f\circ g$ equals that of $g\circ f$ and if they take the same values on it. This is what $f\circ g=g\circ f$ is supposed to mean.
Similarly we will say that ``$f\circ g= h\circ i$, domain included'' if the domain of definition of $f\circ g$ equals that of $h\circ i$ and if they take the same values on it.
We will say that ``$f$ and $g$ commute on $U$'' if $U$ is included in the domain of $f\circ g$ and in the domain of $g\circ f$ and if both compositions take the same value on $U$. There is another natural notion that we will not use: ``$f$ and $g$ commute, wherever defined'', if $f\circ g(z)= g\circ f(z)$ for all $z$ for which both hand sides are defined.

We will denote the domain of definition of a map $f$ by $\on{dom}(f)$.

\section{The construction}

\subsection{Working in the universal cover of $\C^*$}

The map $f$ will be the limit of a sequence of compositions $f_n = g_n \circ \cdots \circ g_1$ of maps $g_i$ fixing the origin. However we prefer to work at the level of lifted coordinates: let 
\[\euh(z)=\exp(2i\pi z).
\]
The map $\euh$ is a universal cover from $\C$ to $\C^*$, with deck transformations $(\Z,+)$ and induces an isomorphism between the cylinder $\C/\Z$ and $\C^*$. This map can be used to extend the Riemann surface $\C/\Z$ at its upper end (we will not use the lower end in this article): denote $+i\infty$ the added point, corresponding to $0\in\C\subset\C^*$.

A holomorphic map $F$ that commutes with $T_1$, domain included, and such that
\[F(z) = z+t+o(1)
\]
for some $t\in\R$ as $\Im(z)$ tends\footnote{The filter understood by ``$\Im(z)\tend +\infty$'' is the one generated by all the half planes $\H_h$ for $h\in\R$. We also imply that the domain of $F$ contains such a half-plane.} to $+\infty$, ``projects'' to a holomorphic map $f$ with an erasable singularity at the origin, i.e. there exists $f$ holomorphic defined on $\{0\}\cup \euh(\text{domain of }f)$, with $\euh\circ F = f\circ \euh$, $f(0)=0$, $f'(0)=e^{2\pi i t}\neq 0$, and $f(z)=0 \implies z=0$.

Conversely, a holomorphic map $f$ defined in an open subset $O$ of $\C$ containing the origin, with $f(0)=0$, $f'(0)\neq 0$, and $f(z)=0 \implies z=0$, will have a (non unique) lift $F: \euh^{-1}(O)\to\C$ that commutes with $T_1$, domain included, with $E\circ F = f\circ E$ and $F(z)-z$ has a limit when $\Im(z)\to+\infty$.

\subsection{Presentation of the actors}\label{subsec:construction}

Let $T_v$ denote the translation by $v$:
\[T_v(z)=z+v.
\]
Let 
\[\H_h=\setof{z\in\C}{\Im(z)>h},\qquad\H=\H_0\]
that we will call \emph{upper half planes}.
Let $\one$ be the constant vector field in $\C$ of expression $\one(z)=1$. 

Assume that maps $\cal R_k$ are given for $k=1$, \ldots, $n$, satisfying the following conditions:
\begin{itemize}
  \item $\cal R_k$ is entire
  \item there exists $q_k\in\N^*$ such that $\cal R_k\circ T_1 = T_{q_k} \circ \cal R_k$
  \item $\cal R_k(z)$ is injective on some upper half-plane. Given the above condition, it is equivalent to ``$\ds \cal R_k(z) = q_k z + c + o(1)$ as $\Im(z)\tend+\infty$''
  \item $\forall z\in\C$, $\cal R_k'(z)\neq 0$
\end{itemize}
We call these maps \emph{renormalizations}. They will be inductively defined later in this article.

Let
\begin{itemize}
\item $X_n = \cal R_1^* \cdots \cal R_n^* \one$, i.e.\ the pull back of the vector field $\one$ by the composition $\cal R_n \circ \cdots \circ \cal R_1$
\item $G_n$, the time $1$ of the flow of the vector field $X_n$
\item $F_n = G_n \circ \cdots \circ G_1$
\end{itemize}
Then $X_n$ is a non-vanishing entire vector field. The maps $G_n$ and $F_n$ are holomorphic but not necessarily defined on all of $\C$.

These objects ``project'' to $\C/\Z$ in the following sense:
\begin{itemize}
\item $T_1^*X_n = X_n$, i.e.\ $X_n$ is $T_1$-invariant
\item $G_n$ and $F_n$ commute with $T_1$, domain included
\item recall that $\cal R_n \circ T_1 = T_{q_n} \circ \cal R_n$
\end{itemize}
In particular, they can be thought of objects living on $\C/\Z$.

Also, they are defined in a neighborhood\footnote{i.e.\ on a set containing a half plane $\H_h$ for some $h$ that may depend on $n$} of $+i\infty$, and we have the following expansions as $\Im(z)\tend+\infty$:
\begin{itemize}
\item $\ds \cal R_n(z) = q_n z + c + o(1)$
\item $\ds X_n(z) = \frac{1}{q_1 \ldots q_n} + o(1)$
\item $\ds G_n(z) = z + \frac{1}{q_1\ldots q_n} + o(1)$
\item $\ds F_n(z) = z + \sum_{k=1}^n \frac{1}{q_1\ldots q_k} + o(1)$
\end{itemize}

Since each $G_n$, and thus $F_n$ commute with $T_1$, domains included, the map $F_n$ projects by $\euh$ to a map $f_n$, i.e. there exists a holomorphic map $f_n:\{0\}\cup \euh(\on{dom}(F_n)) \to \C$ such that $\euh\circ F_n=f_n\circ \euh$, domain included, and $f_n$ fixes the origin.

Let $\Dom$ be the set of points that have a neighborhood on which all $F_n$ are defined and converge uniformly. It is an open set. It might be empty. Let $F$ be the map defined on $\Dom$ as the limit of the $F_n$. Then $F$ is holomorphic and commutes with $T_1$. The following lemma is elementary so we will not explain its proof: it is stated for reference.
\begin{lemma}\label{lem:fEEF}
 Assume that $\Dom$ contains an upper half plane. Let $\dom=\euh(\Dom)\cup\{0\}$. Then $\dom$ is an open subset of $\C$ and $f_n$ converges locally uniformly on $\dom$ to a holomorphic map $f$ defined on $\dom$ and satisfying $f\circ \euh=\euh\circ F$, domain included.
\end{lemma}

The set of domains containing an upper half plane and on which $\cal R_n$ is a bijection to an upper half plane or to $\C$, is non-empty and has a greatest element $V_n$ for inclusion. Let $U_0=\C$ and for $n\geq 1$ let $U_n$ be the set of points $z\in V_1$ such that for all $k<n$, $\cal R_k\circ\cdots\circ\cal R_1(z) \in V_{k+1}$. Then $(U_n)_{n\in\N}$ forms a decreasing sequence for inclusion, of $T_1$-invariant simply connected open subsets $U_n$ of $\C$ containing upper half planes, on which $\cal R_n\circ\cdots\circ\cal R_1$ is an isomorphism to an upper half plane or to $\C$ and on which $\cal R_k\circ\cdots\circ\cal R_1$ is injective for all $k\leq n$. 
Let $\wt U_n=\euh(U_n)\cup\{0\}$. It is a connected and simply connected open subset of $\C$. 
Note that each $f_n$ is conjugated to a rotation in a neighborhood of the origin. Indeed, let $\theta_n=\sum_{k=1}^{n} 1/(q_1\cdots q_k)$.
The map $z\mapsto \frac{1}{q_n}\cal R_n\circ\cdots\circ\cal R_1(z)$ is a conjugacy from $F_n$ on $U_n$ to $T_{\theta_n}$ on an upper half plane or on $\C$, and commutes with $T_1$ on $U_n$, domain included.
Therefore it projects by $\euh$ to a conjugacy on $\wt U_n$ to the rotation of angle $2\pi\theta_n$ on a disk or on $\C$.

We will use the following lemma, applied to well chosen subsets $\wt D_n$ of $\wt U_n$:
\begin{lemma}\label{lem:ee}
Assume $f_n$ is a sequence of holomorphic maps defined on open subsets of $\C$ containing the origin. Assume $f_n$ fixes $0$ and that $f_n'(0)=e^{i2\pi\theta_n}$ with $\theta_n\in\R$ and $\theta_n\tend \theta\in\R$ as $n\to+\infty$. Assume $\wt D_n$ is a sequence of simply connected open subsets of the domain of $f_n$ and containing $0$, such that $f_n$ is analytically conjugated on $\wt D_n$ to the rotation of angle $2\pi\theta_n$ on a disk or on $\C$.
Assume that there exists a limit $\wt D$ to the sequence $\wt D_n$ in the sense of Caratheodory. Then $f_n$ tends on $\wt D$ locally uniformly to a map $f$ that is analytically conjugated on $\wt D$ to the rotation of angle $\theta$ on a disk or on $\C$.
\end{lemma}
\proof Let $R_{\alpha}(z)=e^{i2\pi\alpha}z$.
The conjugacy $\phi_n$ on $\wt D_n$ can always be chosen so that $\phi_n'(0)=1$. Then a theorem of Caratheodory says that $\phi_n$ converges locally uniformly on $\wt D$ to the unique conformal map $\phi$ from $\wt D$ to a disk or $\C$ such that $\phi(0)=0$ and $\phi'(0)=1$. The reciprocal $\phi_n^{-1}$ is locally uniformly convergent to $\phi^{-1}$. From $f_n= \phi_n^{-1} \circ R_{\theta_n} \circ \phi_n$ we deduce that $f_n$ tends to $\phi^{-1}\circ R_{\theta} \circ \phi$ locally uniformly.
\endproof

Depending on authors, the definition of Siegel disks mays differ slightly. For us, it is the maximal domain of conjugacy to an irrational rotation. 
In the lemma above, if $\theta$ is irrational, then $\wt D$ is an invariant subset of the Siegel disk of $f$: could the Siegel disk of $f$ be bigger?
It necessarily requires the boundary of $\wt D$ to be an analytic Jordan curve\footnote{this claim uses the fact that $\theta$ is irrational}, compactly contained in $\dom$. In fact it is equivalent. So if we get a boundary that is not a Jordan curve or not an analytic one, then we know that $\wt D$ equals the Siegel disk of $f$. (This implies that $\wt U_n$ has the same Caratheodory limit as $\wt D_n$.)

How to build examples with $\wt D\subset\subset \dom$? First we will ensure that $\Dom$ contains $\H_{-h}$ for some $h>0$. Second we will ensure that $\wt D$ is contained in $\H$.

\subsection{Ensuring convergence on big domains}\label{subsec:relcomp}

Here we give a slightly informal description, and a rigorous lemma that will be used in the more formal construction described in section~\ref{subsec:core}.

It will be handy to use a notation for the set of maps satisfying the conditions we required on $\cal R_n$ in the previous section. So let $\cal H_1$ denote the set of maps $\cal R$ such that:
\begin{itemize}
  \item $\cal R$ is entire
  \item $\cal R\circ T_1 = T_1 \circ \cal R$
  \item $\cal R(z)$ is injective on some upper half-plane
  \item $\forall z\in\C$, $\cal R'(z)\neq 0$
\end{itemize}
So that $\cal R_k\in q_k\cal H_1$, i.e. $\cal R_k$ is the product of a member of $\cal H_1$ with the scalar $q_k$.

\bigskip

For this, we will inductively define the maps $\cal R_n$ using the following (elementary) lemma:
\begin{lemma}\label{lem:ff}
  Fix $n\in \N$. Assume we are given $\cal R_k\in q_k\cal H_1$ for $k=1$, \ldots, $n-1$ and some $B_n\in \cal H_1$. Let $\cal R_n=qB_n$ for some $q\in\N^*$. The corresponding map $f_n$ depends on $q$ but $f_{n-1}$ does not. Then, as $q\tend+\infty$, $f_n$ tends to $f_{n-1}$ in the following sense: every compact subset $K$ of $\on{dom}(f_{n-1})$ is eventually contained in $\on{dom}(f_{n})$ and $f_{n} \tend f_{n-1}$ uniformly on $K$. Note however, that the domain $U_n$ is independent of $q$.
\end{lemma}

\proof Indeed, $\cal R_n^*\one = (qB_n)^*\one = \frac{1}{q}B_n^*\one$.
So $X_n$ is the product of the scalar $\frac{1}{q}\times$ with an entire vector field that does not depend on $q$. The map $G_n$ is the time $1/q$ of the flow of the latter vector field. The rest follows easily.
\endproof

Fix some number $h>0$. We will inductively define the sequence $\cal R_n$ by first choosing $B_1$, then $q_1$, then $B_2$, then $q_2$, etc\ldots\ We will be careful to choose $q_k$ big enough so that the domain of definition of $f_n$ contains the upper half plane $\H_{-h-\frac{1}{n}}$ and $|f_{n-1}-f_{n}|<1/2^n$ on $\H_{-h}$. This is always possible by the lemma above. We can also choose $q_n$ to ensure that the sum $\theta$ is convergent\footnote{It already follows from the convergence of $f_n$, though.} and has an irrational value: a sufficient condition for this is that $q_n$ tends to $+\infty$. Then the half plane $\H_{-h}$ is contained in the domain of definition $\Dom$ of $F$.

It shall be stressed out that in the construction, the sequence $(B_n)$ will not be fixed independently of the sequence $(q_n)$: the choice of $B_{n+1}$ will depend on the choice of all the previous objects $B_1$, $q_1$, $B_2$, $q_2$, \ldots~$B_{n}$, $q_n$.

\subsection{Ensuring compactly contained bad boundaries}\label{subsec:compact}

As in the previous section we give here a slightly informal description, and a rigorous lemma that will be used in the more formal construction described in section~\ref{subsec:core}.

The following claim is where we use \textbf{Runge's theorem}:
\begin{lemma}\label{lem:runge}
 Assume $\phi$ is a conformal bijection from $\H/\Z$ to a domain 
$W\subset\C/\Z$ that extends to a conformal bijection fixing the upper end.
Then there exists a sequence of holomorphic maps $\psi_n:\C/\Z \to \C/\Z$ such that:
  \begin{itemize}
    \item $\psi'_n$ does not vanish,
    \item $\psi_n(z)-z$ has a limit when $\Im(z)$ tends to $+\infty$,
    \item as $n$ tends to $+\infty$, $\psi_n$ tends uniformly to $\phi^{-1}$ on every compact subset of $W$ and in a neighborhood of the upper end of $\C/\Z$.
  \end{itemize}
\end{lemma}
  \proof
  Conjugate the situation by $\euh$: let $\wt W = \euh(W)\cap\{0\}$ and $\Psi: \wt W \to \D$ be the conjugate of $\phi^{-1}$, extended by setting $\Psi(0)=0$, where $\D$ denotes the unit disk. It is a conformal isomorphism. It can therefore be written as $\Psi(z)=ze^{v(z)}$ for some holomorphic map $v: \wt W \to \C$. As a conformal isomorphism its derivative does not vanish, hence $1+z v'(z)\neq 0$. Thus there exists a map $u: \wt W \to \C$ such that $1+z v'(z) = e^{u(z)}$. Since $e^{u(0)}=1$, we can choose $u$ such that $u(0)=0$. So $u(z)=zw(z)$ for some holomorphic $w: \wt W \to \C$.
Now apply Runge's theorem: there exists a sequence of entire maps $w_n$ (even polynomials) that converge locally uniformly to $w$ on $\wt W$. There exists an  entire map $v_n$ (unique) taking the same value as $v$ at the origin and such that $1+zv'_n=e^{zw_n}$. Indeed $(e^{zw_n}-1)/z$ has a removable singularity at the origin, so it extends to an entire map, of which $v_n$ is just the appropriate primitive. The maps $\psi_n(z)=z+v_n(e^{z})$ satisfy the required conditions.
\endproof

Now we describe informally how we can use this lemma to choose the functions $B_n$ to get interesting examples. A more formal, but more specific, description is done in section~\ref{subsec:core}.

Let $D_0=\H/\Z$.
Let $W_1$ be any open simply connected strict subset of $\H/\Z\cup\{+i\infty\}$ containing the upper end $+i\infty$. Let $\phi$ be any conformal map from $\H/\Z\cup\{+i\infty\}$ to $W_1$ that fixes the upper end.
Let $\epsilon>0$ and let the Jordan curve $J_1$ be the image by $\phi$ of the Jordan curve of equation ``$\Im(z)=\epsilon$'' in $\H/\Z$.
If the boundary of $W_1$ is a convoluted Jordan curve (or something more complicated than a Jordan curve) and $\epsilon$ is small then $J_1$ will also be convoluted.

Apply lemma~\ref{lem:runge} to $W=W_1-\{+i\infty\}$.
By the uniform convergence of $\psi_n$ to $\phi^{-1}$ on compact subsets of $W$ and on a neighborhood of the upper end of $\C/\Z$, there is for $n$ big enough a branch $\phi_n$ of $\psi_n^{-1}$ defined for $\Im(z)>\epsilon/2$ and mapping the upper end to itself. For $n$ big, the Jordan curve $J'_1=\phi_n(\R+i\epsilon)$, image by $\phi_n$ of the Jordan curve of equation ``$\Im(z)=\epsilon$'' in $\H/\Z$, will be close to $J_1$ so it will also be convoluted.

So we will choose $B_1=\psi_n(z)-i\epsilon$ for some big value of $n$. Then we will chose $q_1$ as in section~\ref{subsec:compact}. Let $D_1$ be the domain bounded by $J'_1$ and containing the upper end.
Among the integral lines of the vector field $X_1$, there are the images of the horizontals by $\phi_n$, which loop. This includes $J'_1$ and also a set of curves foliating the domain $D_1$. The map $F_1=G_1$, seen as acting on $\C/\Z$, is conjugate to a finite order\footnote{Recall that these maps are not entire. Some iterate will be the identity on the component of the domain for the iterate that contains $D_1$.} rotation on $D_1$. The renormalization $\cal R_1$ induces a conjugacy of the restriction of $F_1$ to $D_1$ to the translation $T_1$ on the half cylinder $\H/q_1\Z$.

The next step in the construction consists in choosing a new simply connected domain $W_2$, contained in $D_1$, containing the upper end, whose boundary is even more convoluted, yet invariant by $F_1$, and very close to $J'_1$. By the conjugacy $\cal R_1$ on $D_1$ it amounts to choosing a domain in $\H/q_1\Z\cup\{+i\infty\}$ that is invariant by $T_1$, i.e. to choosing a simply connected domain $W_2'$ in $\H/\Z\cup\{+i\infty\}$ (by applying the cover $z+q_1\Z\in\H/q_1\Z \mapsto z+\Z \in \H/\Z$), and whose boundary is very close to $\R$.

Let $\phi : \H/\Z\cup\{+i\infty\} \to W_2'$ be a conformal map fixing the upper end.
We choose a new value of $\epsilon$, which gives a new Jordan curve $J_2=D_1\cap \cal R_1^{-1}(\phi(\R+i\epsilon))$ close to the boundary of $W_2$. We apply lemma~\ref{lem:runge} again and get a new sequence of entire maps $\psi_n$ tending to $\phi^{-1}$ on $W_2'$. We set $B_2=\psi_n(z)-i\epsilon$ for $n$ big enough so that the curve $J'_2=D_1\cap \cal R_1^{-1}(\phi_n(\R+i\epsilon))$ is close to $J_2$.
We choose $q_2$ as in section~\ref{subsec:compact}.
The domain $D_2$ bounded by $J'_2$ is then invariant by $G_1$ and by the flow of the vector field $X_2$. The maps $G_1$ and $G_2$, will be thus conjugated by $\cal R_2 \circ \cal R_1$ on $D_2$ to the translations by respectively $q_2$ and $1$ on $\H/q_1q_2\Z$.

And so on\ldots\ We get a decreasing sequence of domains $D_n \subset\H/\Z$ on which $\cal R_n \circ \cdots\circ \cal R_1$ is a conformal bijection to $\H/(q_1\cdots q_n)\Z$ and conjugating all the maps $G_k$, for $k\leq n$, to integer translations, and thus the map $F_n$ too. Their boundaries form a sequence of Jordan curves extremely convoluted and close to each other. The construction has enough flexibility to allow for $D=\bigcap_n D_n$ to have a non locally connected boundary, as we will prove in the following sections (basically by the same method as Handel in \cite{handel}).

\subsection{Pseudo circles}\label{subsec:pseudo}

Pseudo circles are amazing topological objects. The interested\footnote{Let us mention the following facts: There exists a circularly chainable and hereditarily indecomposable planar continuum. For the definition of all these terms, the reader is referred to the introduction of \cite{rogers}. It was proved to be unique up to homeomorphism (but it is not homogeneous). Such an object is called a pseudo circle.} reader may look at \cite{bing} (Example 2 p.48), \cite{fearnley}, \cite{rogers2}, \cite{KK}. It was chosen by Handel in \cite{handel} as a good example of extreme pathology that a minimal set can have\, even for an area preserving smooth diffeomorphism of the sphere. Herman did the same choice in \cite{herman} for his construction. So we will continue this tradition.

For our purpose, it will be enough to know the following:

\begin{itemize}
  \item A pseudo circle is compact, connected, separates the plane into two components, but it is not locally connected thus it is not a Jordan curve.
  \item It is equal to the boundary of both components (in particular it has empty interior).
  \item Any set $K$ obtained by the following procedure is a pseudo circle.
\end{itemize}

\textsc{Procedure}:
A \emph{circular chain} in a topological space $X$ is a sequence $e=(e_i)$ of open subsets of $X$, called \emph{links}, indexed by $\Z/m\Z$ for some $m\geq 4$, such that $e_i\cap e_j\neq\emptyset$ if and only if $i,j$ are adjacent or equal.\footnote{This is probably inspired from \v{C}ech cohomology: we want the \emph{nerve} of the open cover to be a circle. So for chains of length $3$ one would also require $e_{\overline{0}} \cap e_{\overline{1}} \cap e_{\overline{2}} = \emptyset$.}
For each $n\geq 1$, let $Q_n$ be a circular chain in the plane\footnote{it does not matter whether the $e_i$ are connected or not} such that:
\begin{itemize}
  \item the diameters of the links of $Q_n$ are finite for all $n$ and their supremum for a given $n$ tends to $0$ as $n\to+\infty$
  \item $Q_{n+1}$ is a refinement of $Q_n$: each link of the former is contained in a link of the latter
  \item $Q_{n+1}$ is \emph{crookedly embedded} in $Q_n$, which means the following: assume that $Q_n=(d_i)$ has $m$ links, and that $Q_{n+1}=(e_i)$ has $m'$ links. It will be convenient to let $i\in \Z$ and let $d_i$ denote $d_{\bar{i}}$ where $\bar{i}$ is the class of $i$ modulo $m$. Similarly, for $i\in\Z$ let $e_i$ denote $e_{\bar{i}}$ for the residue class modulo $m'$. To be \emph{crookedly embedded} means that there exists a map $f:\Z\to\Z$ such that 
  \begin{enumerate}
    \item\label{item:cr1} $\forall i\in\Z$, the closure of $e_i$ is contained in $d_{f(i)}$
    \item\label{item:cr4} $\forall i\in\Z$, $f(i+1)-f(i)\in\{-1,0,1\}$ (\footnote{By condition \ref{item:cr1} it is already the case modulo $m$ because $Q_n$ and $Q_{n+1}$ are circular chains.})
    \item\label{item:cr2} $\forall i\in\Z$, $f(i+m')=f(i)+m$ (\footnote{It translates the idea of $Q_{n+1}$ having winding number $1$ in $Q_n$.})
    \item\label{item:cr3} for all $i,j$ with $f(i)+2< f(j)<f(i)+m$, there exists $i'$ and $j'$ with $f(i')=f(j)-1$, $f(j')=f(i)+1$ and either $i<i'<j'<j$ or $i>i'>j'>j$. 
  \end{enumerate}
\end{itemize}
Then define $K=\bigcap_{n\geq 1} Q_n$: it is a pseudo circle.\hfill\textsc{end of description}

\medskip

Let us say that a Jordan curve, parameterized by $\iota(t)$, $t\in\R/\Z$, is \emph{stably crooked} in a circular chain $Q_n=(e_i)$ with $i\in\Z/m\Z$, whenever:
\begin{itemize}
  \item the curve is contained in the union of the links of the chain,
  \item it has winding number $1$ in $Q_n$, oriented by $\iota$, (\footnote{The winding number of an oriented closed curve $\iota$ in a circular chain $Q=(d_j)$ of length $m\geq 4$ could be defined in terms of first \v{C}ech cohomology groups with coefficients in $\Z$ by a diagram looking like this: $\Z \approx \check{H}^1(Q,\Z) \mapsto \check{H}^1(\bigcup_\jbar d_\jbar,\Z) \mapsto \check{H}^1(S^1,\Z) \approx \Z$, but this is sophisticated. 
  It amounts to doing the following: cut the curve into a finite number of small consecutive pieces $e_\ibar=\iota([t_\ibar,t_{{\overline{\imath+1}}}])$ for some sequence $t_\ibar$ indexed by $\Z/p\Z$ for some $p>0$. If they are small enough, they will all be contained in some $d_\jbar$. There is at most two possible values of $\jbar$. Choose anyone and call it $g(\ibar)$. Since $Q_n$ is a chain, $g({\overline{\imath+1}})$ and $g(\ibar)$ must differ of at most one in $\Z/m\Z$. There exists a lift $f:\Z\to\Z$ of $g$ such that $f(i+1)$ and $f(i)$ differ of at most one for all $i$. Then there exists $k\in\Z$ such that $f(i+p)=f(i)+km$ for all $i$. This $k$ is the winding number. It is independent of all choices and of the oriented parameterization $\iota$, and invariant by small perturbations of $\iota$.})
  \item $\forall t_1,t_2\in\R$ with $t_1<t_2<t_1+1$ and $\forall u,v\in\Z$, if $\iota(t_1)\in \overline{e}_u$ and $\iota(t_2)\in \overline{e}_v$ and if $u+2< v < u+m$ then there exists $t'_1$ and $t'_2$ with either $t_1<t'_1<t'_2<t_2$ or $t_1>t'_1>t'_2>t_2$ such that $\iota(t'_1)\in e_{v-1}$ and $\iota(t'_2)\in e_{u+1}$.
\end{itemize}
It is independent of the choice of parameterization of the Jordan curve.
The subtlety of taking closures for $\overline{e}_u$ and $\overline{e}_v$ but not $e_{u+1}$ and $e_{v-1}$ is to ensure the following: \emph{if a Jordan curve is stably crooked in $Q_n$ then any nearby Jordan curve is}. By nearby we mean a Jordan curve which is parameterized by a map that is close to $\iota$.

\subsection{The core of the construction}\label{subsec:core}

\begin{figure}
\scalebox{0.5}{\includegraphics{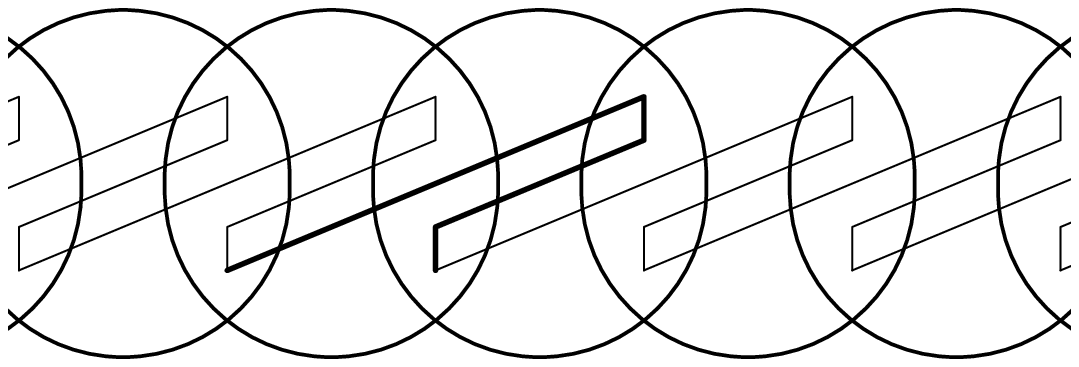}}\\
\scalebox{0.5}{\includegraphics{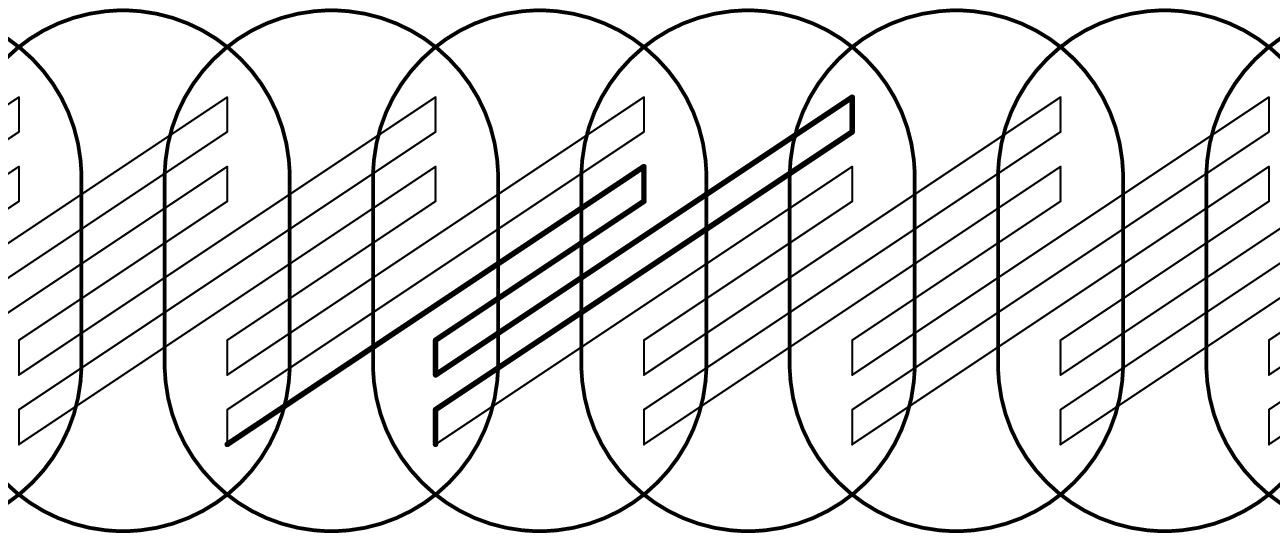}}\\
\scalebox{0.5}{\includegraphics{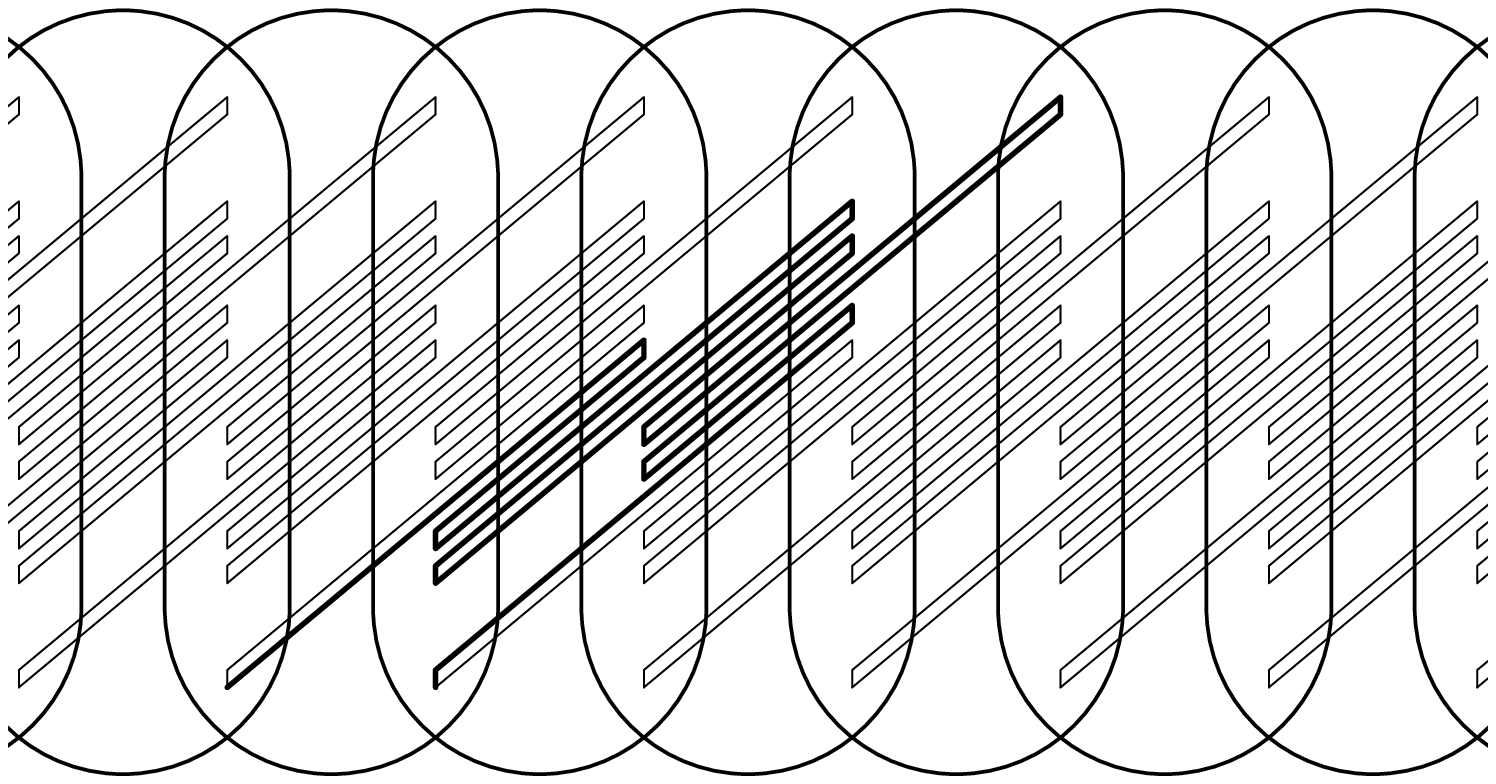}}
\caption{Examples of $T_1$-invariant stably crooked Jordan curves in circular chains (unrolled) of length $5$, $6$ and $7$. To make the figure more readable, we replaced the rectangles of the chain $C_m$ by smoother regions.}
\label{fig:crook1}
\end{figure}

Following Herman (with different conventions), for $m\geq 4$ consider the circular chain $C_m$ of length $m$ consisting of the rectangular links $e_i=\big]i-1/4, i+5/4\big[ \times \big]0, 1\big[$ in the cylinder $\C/m\Z$.
Then
\begin{itemize}
  %
  %
  \item \emph{There exists a Jordan curve in $\C/m\Z$ that is stably crooked in $C_m$ and invariant by $T_1$.} The proof is in \cite{handel} but let us mention that it is quite entertaining to work it out on one's own. Figure \ref{fig:crook1} gives examples of solutions for $m=5$, $6$ and $7$.
  \item Fix such a Jordan curve $J$. Consider a conformal map $\phi$ from $\H/m\Z\cup\{+i\infty\}$ to the component of $\C/m\Z\cup\{+i\infty\} - J$ that contains the upper end, and that fixes the upper end. There exists such maps.\footnote{Apply the Riemann mapping theorem to the image of the situation by $\euh$.} \emph{There exists $\epsilon$ such that the Jordan curve $J_{\epsilon}$, image by $\phi$ of the curve of equation ``$\Im(z)=\epsilon$'', is stably crooked in $C_m$.} It follows from the definition of stably crooked and from the following theorem of Caratheodory: $\phi$ extends continuously to the closure of $\H/m\Z$ and the extension is on $\R/m\Z$ a parameterization of the Jordan curve $J$.
  \item Let $W$ be the preimage by the natural projection $\C\to \C/m\Z$ of the image of $\phi$. Then $W$ is connected and $\phi$ lifts to a conformal bijection, also denoted by $\phi$, from $\H$ to $W$. Since $W$ is $T_1$-invariant, $\phi$ commutes with $T_1$, domain included. So $\phi$ induces an isomorphism from $\H/\Z$ to $W/\Z$ that we will also denote by $\phi$.
  \item Now apply lemma~\ref{lem:runge} to $\phi$ and $W/\Z$. We get a sequence of holomorphic maps $\psi_n:\C/\Z \to \C/\Z$ such that:
  \begin{itemize}
    \item $\psi'_n$ does not vanish,
    \item $\psi_n(z)-z$ has a limit when $\Im(z)$ tends to $+\infty$,
    \item as $n$ tends to $+\infty$, $\psi_n$ tends uniformly to $\phi^{-1}$ on every compact subset of $W$ and in a neighborhood of the upper end of $\C/\Z$.
  \end{itemize}
  \item  By this uniform convergence, there is for $n$ big enough a branch $\phi_n$ of $\psi_n^{-1}$ defined for $\Im(z)>\epsilon/2$ and mapping the upper end to itself.
  \emph{There exists $M\geq 4$ and $N\in\N$ such that for all $m'\geq M$ and $n\geq N$, the image of the chain $C_{m'}$ by $z\mapsto\phi_n\big(\frac{m}{m'}z+i\epsilon\big)$ is a circular chain that is crookedly embedded in $C_m$. }
  \par\smallskip\noindent\textit{Proof.} Denote $(e_\ibar)$ this circular chain and $C_m=(d_\jbar)$, with $\ibar\in\Z/m'\Z$ and $\jbar\in\Z/m\Z$.
   First there is $M, N$ such that for $m'\geq M$ and $n\geq N$, for all link $e_\ibar$, the set of $\jbar$ such that ``$\overline{e_\ibar}\subset d_\jbar$'' is not empty, for otherwise, taking subsequences, there would be a point of $J_{\epsilon}$ that is a limit of points not contained in any link, but links are open and $J_{\epsilon}$ is contained in their union, which leads to a contradiction. Since $C_m$ is a chain, this set of $\jbar$ is at most two consecutive integers modulo $m$. Let $g(\ibar)$ be any of them (the ``smallest'' for instance). Because $(e_\ibar)$ and $(d_\jbar)$ are circular chains, we have necessarily $g(\overline{\imath+1})-g(\ibar)\in\{\overline{-1},\overline{0},\overline{1}\}$. There is thus a unique lift of $f$ from $\Z$ to $\Z$ that satisfies $f(i+1)-f(i)\in\{-1,0,1\}$ and $f(0)\in\{0,\ldots, m-1\}$. This lift necessarily satisfies $f(i+m')=f(i)+km$ for some $k$ independent of $i$. This $k$ is necessarily equal to the winding number in $C_m$ of the curve $\phi_n(x+i\epsilon)$ that runs along the lower edges of the ``rectangles'' defining $(e_\ibar)$, which is equal to the winding number of $J_\epsilon$ in $C_m$, that is $1$, so we get condition \eqref{item:cr2}.
If condition \eqref{item:cr3} were not satisfied for all $m',n$ big enough, then taking subsequences would contradict the fact that $J_\epsilon$ is stably crooked in $C_m$.
  \qed
\end{itemize}

Let us sum up what we will use of the above analysis. Setting $\psi(z)=\psi_N(z)-i\epsilon$ for the value of $N$ mentioned above, we get:
\begin{lemma}\label{lem:psi}
  For all $m\geq 4$ there exists $M\geq 4$ and a  holomorphic map $\psi$ such that 
  \begin{itemize} 
    \item $\psi$ is entire,
    \item $\psi'$ does not vanish,
    \item $\psi\circ T_{1} = T_{1} \circ \psi$,
    \item $\psi(z)-z$ has a limit when $\Im(z)\to+\infty$,
    \item $\psi$ has an inverse branch $\phi$ defined for $\Im(z)>-\epsilon/2$ for some $\epsilon>0$, with $\phi(z)-z$ having a limit as $\Im(z)\to+\infty$, and $\phi\circ T_1=T_1\circ\phi$,
    \item for all $m'\geq M$, the image of the circular chain $C_{m'}$ by $z\mapsto\phi\big(\frac{m}{m'} z\big)$ is a circular chain crookedly embedded in $C_m$.
  \end{itemize}
\end{lemma}

\subsection{Putting it all together}

We can now apply this to build our example. Fix $h>0$. Let $B_n \in\cal H_1$, $q_n\in\N^*$ be defined by induction as follows:

\smallskip

Let $B_1=\on{id}$ on $\C$, and $q_1$ be any integer $\geq 5$.

\smallskip

Assume that $B_1$, $q_1$, \ldots\ $B_{n}$, $q_{n}$ have been fixed such that the following holds for all $k\leq n$:
\begin{enumerate}
  \item\label{item:1} As in section~\ref{subsec:construction}, let $\cal R_k=q_k B_k$, let $V_k$ be biggest domain containing an upper half plane and on which $\cal R_k$ is a bijection to an upper half plane, let $U_0=\C$ and let $U_n$ be the set of points $z\in V_1$ such that for all $k<n$, $\cal R_k\circ\cdots\circ\cal R_1(z) \in V_{k+1}$. We already saw that $\cal R_k\circ\cdots\circ\cal R_1$ are bijections from $U_k$ to upper half planes. We require these half-planes to contain $\H_{-\epsilon}$ for some $\epsilon>0$ which may depend on $k$. 
  \item\label{item:2} Let $Q_k$ be the circular chain defined as the preimage in $\C/\Z$ by the restriction of $\cal R_k\circ\cdots\circ\cal R_1$ to $U_k/\Z$ of the chain $C_m\subset\H/m\Z$ with $m=q_1\ldots q_{k}$. If $k>1$, we require $Q_k$ to be crookedly embedded in $Q_{k-1}$.
  \item\label{item:3} The supremum of the diameters of the links of $Q_k$ is $\leq 1/k$. 
  \item\label{item:4} Let $F_k$ be defined as in section~\ref{subsec:construction}. We require the domain of definition of $F_k$ to contain $\H_{-h-1/k}$.
  \item\label{item:5} If $k>1$, the supremum on $\H_{-h}$ of $\big|F_k-F_{k-1}\big|$ is $\leq 1/2^k$.
  \item\label{item:6} $q_k\geq k$.
\end{enumerate}

\smallskip

This set of conditions is satisfied for $n=1$: for \eqref{item:1} we have $\cal R_1(z)=5z$, $\C=V_1=U_1=\cal R_1(U_1)$; \eqref{item:2} is empty; for \eqref{item:3} the diameter of the links $C_m$ is $\sqrt{13}/2$ thus the diameter of the links of $Q_1$ is $\sqrt{13}/2q_1$ and here $m=q_1\geq 5$; \eqref{item:4} $F_1=T_{1/q_1}$ and is defined on $\C$; \eqref{item:5} is empty; \eqref{item:6} $q_1\geq 5\geq 1$.

\smallskip

If it is satisfied for some $n$, then let us explain why it is possible to choose $(B_{n+1},q_{n+1})$ such that it is satisfied for $n+1$:
\begin{itemize}
\item First, we define $B_{n+1}$. Let $m=q_1 \cdots q_n$.
Let $\psi$, $\phi$, $M$ be given by lemma~\ref{lem:psi}.
Let $B_{n+1}(z)=\psi(z)$ and for any $q_{n+1}\in\N^*$ set $m'=q_{n+1}m$.
We obtain at once \eqref{item:1}, and \eqref{item:2} follows as soon as $q_{n+1}$ is big enough so that $m'\geq M$.
\item That \eqref{item:3} holds when $q_{n+1}$ is big enough follows merely from the continuity of $\phi$.
\item That \eqref{item:4} and \eqref{item:5} hold when $q_{n+1}$ is big enough follows from section~\ref{subsec:relcomp}.
\item Point \eqref{item:6} needs no comment.
\end{itemize}

Then by \eqref{item:5} the sequence $F_n$ tends uniformly on $\H_{-h}$ to some holomorphic map $F$, so $\H_{-h}\subset\Dom$ in the notations of section~\ref{subsec:construction}. Let $f$ the map associated to $F$ in  lemma~\ref{lem:fEEF}. 
By \eqref{item:6} the number $\theta=\sum_n 1/(q_1\cdots q_n)$ is irrational. By \eqref{item:2} and \eqref{item:3} the set $K=\bigcap_n Q_n$ is a pseudo circle.
Let $D_n$ be the preimage of $\H$ by the restriction of $\cal R_k\circ\cdots\circ\cal R_1$ to $U_k$ and let $\wt D_n=\euh(D_n)\cup\{0\}$. By \eqref{item:1}, $f_n$ is conjugated to a (rational\footnote{Some iterate of $f_n$ is thus the identity on the component of its domain that contains $U_n$. Recall that $f_n$ is not entire: when $f_n$ is iterated, its domain decreases and may disconnect.}) rotation on $\wt D_n$ (see the discussion in section~\ref{subsec:construction}). The sequence $\wt D_n$ is decreasing for inclusion, and tends in the sense of Caratheodory to 
the connected component $\wt D$ of the complement of the pseudo circle $\euh(K)$ that contains the origin. The boundary of $\wt D$ is equal to $E(K)$ (it is a property that pseudo circles share with Jordan curves: they are equal to the boundary of both components of their complement).
By lemma~\ref{lem:ee}, $\wt D$ is contained in the Siegel disk of $f$. But since its boundary is not a Jordan curve (it is a pseudo circle) and the rotation number is irrational, $\wt{D}$ \emph{is} the Siegel disk. It is contained in the unit disk $\D$, because $\wt D_1=\H$ and $\wt D_n$ is decreasing. So the Siegel disk of $f$ is compactly contained in the domain of definition of $f$, because the latter contains $e^{2\pi h}\D$.

\section{Conclusion}

With the above proof, we realize that the map $f$ we constructed is injective, as a limit of injective holomorphic maps. The domain of definition can be taken to contain $R\D$ for $R$ as big as wanted (while the pseudo circle is contained in $\D$). By injectivity and properties of univalent maps, $f$ will be very close to a rotation on, say, $2\D$ when $R$ gets big. The pseudo circle can be chosen to be very close to the unit circle. Alternatively, with a slight modification of the construction it can be chosen to span between distance $\epsilon$ and $1-\epsilon'$.

In \cite{ricardo}, it is mentioned that with Pérez Marco's examples, one can get a whole uncountable group of commuting maps sharing the same Siegel disk. Here the same holds: by taking the $u_n$ small enough we can ensure that for any infinite subset $J$ of $\N$, the infinite composition of the $G_k$ over $k\in J$ will converge on big domains, and they will all leave the Siegel disk of $f$ invariant, so they will have the same Siegel disk since their rotation number remains irrational.

We believe that the new flexibility allowed by Runge's theorem allows to prove the following claims. But this has to be carefully checked.

\begin{claim} Prove that there exists an injective holomorphic map $f$ defined in a simply connected open subset $U$ of $\C$ containing the origin, fixing $0$ and having at $0$ a hedgehog of positive Lebesgue measure compactly contained in $U$.
\end{claim}

\begin{claim}
Prove that there exists an injective holomorphic map $f$ defined in a doubly connected open subset $U$ of $\C$ and a Jordan curve $J$ with positive Lebesgue measure contained in $U$ that is invariant by $f$, and carries an invariant line field.
\end{claim}

We still believe in the following conjecture:

\begin{conjecture}
The boundary of the Siegel disks of all polynomials are Jordan curves.
\end{conjecture}

It seems very likely that the construction requires gigantic values of $q_k$ so that we probably always get very Liouvillian rotation numbers (in particular non Brjuno). We can always arrange so that it is the case.
So Herman's question in \cite{herman} concerning which rotation numbers allow for non Jordan curve boundaries is still quite topical.

\bibliographystyle{alpha}
\bibliography{pseudo}

\end{document}